\documentclass[12pt,twoside,reqno]{amsart}   

\author{Andreas-Stephan Elsenhans}
\address{School of Mathematics and Statistics F07,
              University of Sydney,
              NSW 2006,
              Sydney, Australia}
\email{stephan@maths.usyd.edu.au}
\urladdr{http://www.maths.usyd.edu.au/ut/people?who=AS\_Elsenhans}

\author{J\"org Jahnel}
\address{D\'epartement Mathematik,
              Universit\"at Siegen,
              Walter-Flex-Str.~3,
              D-57068 Siegen, Germany}
\email{jahnel@mathematik.uni-siegen.de}
\urladdr{http://www.uni-math.gwdg.de/jahnel/}


\usepackage{amssymb}
\usepackage{amsmath}
\usepackage{bbm}
\usepackage{mathrsfs}
\usepackage{xypic}

\usepackage{amsthm}
\theoremstyle{plain}
 \newtheorem{theorem}{Theorem}[section]
 \newtheorem{lemma}[theorem]{Lemma}
 \newtheorem{corollary}[theorem]{Corollary}
 \newtheorem{proposition}[theorem]{Proposition}

\theoremstyle{definition}
 \newtheorem{ttt}[theorem]{}
 \newtheorem{nota}[theorem]{Notation}

\theoremstyle{remark}
 \newtheorem{example}[theorem]{Example}
 \newtheorem{exs}[theorem]{Examples}
 \newtheorem{remark}[theorem]{Remark}
 \newtheorem{rems}[theorem]{Remarks}

\newcommand{\calH}{\mathscr{H}}
\newcommand{\calO}{\mathscr{O}}
\newcommand{\calT}{\mathscr{T}}

\newcommand{\bbC}{{\mathbbm C}}
\newcommand{\bbF}{{\mathbbm F}}
\newcommand{\bbN}{{\mathbbm N}}
\newcommand{\bbQ}{{\mathbbm Q}}
\newcommand{\bbR}{{\mathbbm R}}
\newcommand{\bbZ}{{\mathbbm Z}}

\newcommand{\frakm}{\mathfrak{m}}

\newcommand{\id}{\mathop{\textrm{\upshape id}}}
\newcommand{\rk}{\mathop{\textrm{\upshape rk}}}
\newcommand{\Hom}{\mathop{\textrm{\upshape Hom}}\nolimits}
\newcommand{\disc}{\mathop{\textrm{\upshape disc}}}

\newcommand{\cha}{\mathop{\textrm{\upshape char}}}
\newcommand{\Frob}{\mathop{\textrm{\upshape Frob}}\nolimits}
\newcommand{\Gal}{\mathop{\textrm{\upshape Gal}}}
\newcommand{\Pic}{\mathop{\textrm{\upshape Pic}}}
\newcommand{\Sq}{\mathop{\textrm{\upshape Sq}}\nolimits}
\newcommand{\Ann}{\mathop{\textrm{\upshape Ann}}\nolimits}
\newcommand{\Kum}{\mathop{\textrm{\upshape Kum}}\nolimits}

\renewcommand{\pmod}[1]{\nobreak\ifinner\mkern8mu\else\mkern18mu\fi
 (\textrm{\upshape{mod}}\,\,#1)}

\renewcommand{\dim}{\mathop{\textrm{\upshape dim}}\nolimits}
\renewcommand{\deg}{\mathop{\textrm{\upshape deg}}\nolimits}
\renewcommand{\det}{\mathop{\textrm{\upshape det}}\nolimits}
\renewcommand{\ker}{\mathop{\textrm{\upshape ker}}\nolimits}

\newcommand{\N}{\mathop{\textrm{\upshape N\,\!}}}
\newcommand{\tr}{\mathop{\textrm{\upshape tr}}\nolimits}
\newcommand{\Q}{\mathop{\textrm{\upshape Q}}}

\newcommand{\bP}{{\textbf{P}}}
\newcommand{\F}{{\textbf{\upshape F}}}

\newcommand{\et}{{\textrm{\upshape \'et}}}
\newcommand{\tors}{{\textrm{\upshape tors}}}

\renewcommand{\atop}[2]{\genfrac{}{}{0pt}{}{#1}{#2}}

\makeatletter
\def\hsmash{\relax 
  \ifmmode\def\next{\mathpalette\mathhsm@sh}\else\let\next\makehsm@sh
  \fi\next}
\def\makehsm@sh#1{\setbox\z@\hbox{#1}\finhsm@sh}
\def\mathhsm@sh#1#2{\setbox\z@\hbox{$\m@th#1{#2}$}\finhsm@sh}
\def\finhsm@sh{\wd\z@\z@ \box\z@}
\makeatother

\newcounter{abc}
\newenvironment{abc}{\begin{list}{\rm \alph{abc}) }%
{\usecounter{abc} \leftmargin=0.0pt \labelsep=0.0pt %
\listparindent=0.0pt \labelwidth=0.0pt \parsep=\smallskipamount %
\itemsep=0.0pt \topsep=0.0pt \partopsep=\smallskipamount}}{\end{list}}
\newcounter{iii}
\newenvironment{iii}{\begin{list}{\rm \roman{iii}) }%
{\usecounter{iii} \leftmargin=0.0pt \labelsep=0.0pt %
\listparindent=0.0pt \labelwidth=0.0pt \parsep=\smallskipamount%
 \itemsep=0.0pt \topsep=0.0pt \partopsep=\smallskipamount}}{\end{list}}

\def\rightend#1#2{{%
 \leavevmode\nobreak\hskip .5em plus 1fil
 \penalty600 \hskip 0pt plus -1filll
 \vadjust{}\nobreak\hskip 0pt plus 1filll%
 #1\parfillskip=#2\relax \par}}

\def\eop{\ifmmode\rule[-22pt]{0pt}{1pt}\ifinner\tag*{$\square$}\else\eqno{\square}\fi\else\rightend{$\square$}{0pt}\fi}

\title[On the characteristic polynomial of Frobenius]{On the characteristic polynomial of the Frobenius on \'etale cohomology}

\begin{document}

\begin{abstract}
Let~$X$
be a smooth proper variety of even dimension
$d$
over a finite~field. We~establish a restriction on the value
at~$(-1)$
of the characteristic polynomial of the Frobenius on the middle-dimensional \'etale cohomology
of~$X$
with coefficients
in~$\bbQ_l(d/2)$.
\end{abstract}

\maketitle

\section{Introduction}
\label{intr}

Let~$X$
be a smooth proper variety over a finite
field~$\bbF_{\!q}$
of
characteristic~$p > 0$.
Then~the geometric Frobenius
$\Frob$
operates linearly on the
\mbox{$l$-adic}
cohomology vector spaces
$\smash{H^i_\et(X_{\overline\bbF_{\!q}}\!, \bbQ_l(j))}$.
The~characteristic polynomial
of~$\Frob$
has rational coefficients~\cite[Th\'eor\`eme~(1.6)]{De1}.

By~far not every polynomial
$\Phi \in \bbQ[T]$
may occur as the characteristic polynomial of the Frobenius on a smooth proper~variety. The~following conditions were established essentially in the 70s of the last century, in the course of working out A.\,Grothen\-dieck's ideas in Algebraic~Geometry.

\begin{theorem}[{\rm Deligne, Mazur, Ogus}{}]
\label{class}
Let\/~$X$
be a smooth proper variety over a finite
field\/~$\bbF_{\!q}$
of
characteristic\/~$p$.
For\/~$i,j \in \bbZ$,
denote by\/
$$\smash{\Phi_j^{(i)} = T^N + a_1^{(i)} T^{N-1} + \cdots + a_{N-1}^{(i)} T + a_N^{(i)}}$$
the characteristic polynomial
of\/~$\Frob$
on\/~$\smash{H^i_\et(X_{\overline\bbF_{\!q}}\!, \bbQ_l(j))}$.
Then,~for every\/
$r$,
one has\/
$\smash{a_r^{(i)} \in \bbQ}$
and this value is independent
of\/~$l \neq p$.
Moreover,

\begin{abc}
\item
every complex zero
of\/~$\smash{\Phi_j^{(i)}}$
is of absolute
value\/~$q^{i/2 - j}$.
\item
If\/~$i$
is odd then all real zeroes
of\/~$\smash{\Phi_j^{(i)}}$
are of even~multiplicity.
\item
For~every\/
$l \neq p$,
the zeroes
of\/~$\smash{\Phi_j^{(i)}}$
are\/
\mbox{$l$-adic}~units,
i.e.~units in a suitable extension
of\/~$\bbZ_l$.
\item
Put\/~$h_{i-m,m} \!:=\! \dim H^{i-m} (X, \Omega_X^m)$
and define the step-function\/
$G^{(i)} \colon [0,N] \to \bbR$
by
$$
G^{(i)}(t) = \left\{
\begin{array}{ll}
0 & \quad {\rm ~for~} t \leq h_{i,0} \, , \\
n & \quad {\rm ~for~} h_{i,0} + \cdots + h_{i-n+1,n-1} < t \leq h_{i,0} + \cdots + h_{i-n,n} \, .
\end{array}
\right.
$$
Then,~for\/
$r = 1, \ldots, N$,
one has\/
$a_r^{(i)} = 0$~or
$$\nu_q (a_r^{(i)}) \geq \int\limits_0^r [G^{(i)}(t) - j] \,dt \, .$$
Here,
$\nu_q$
is the non-archimedean valuation such
that\/~$\nu_q(q) = 1$.
\end{abc}
\end{theorem}

\begin{rems}
\begin{iii}
\item
We will provide references and proof sketches at the beginning of Section~\ref{sec_proo}.
\item
Assertion~a) immediately implies that
$\smash{\Phi_j^{(i)} \in \bbQ[T]}$
fulfills the functional~equation
$$T^N \Phi(q^{i-2j}/T) = \pm q^{\frac{N}2(i - 2j)} \Phi(T)$$
for~$\smash{N := \dim H^i_\et(X_{\overline\bbF_{\!q}}\!, \bbQ_l(j))}$.
Indeed,~on both sides, there are polynomials with leading
term~$\smash{\pm q^{\frac{N}2(i - 2j)} T^N}$.
They~have the same zeroes as,
with~$z$,
the number
$\smash{\overline{z} = \frac{q^{i-2j}}z}$
is a zero
of~$\Phi$,~too.

Furthermore,~by~b), the plus sign always holds
when~$i$
is~odd. The~statements c) and~d) together show
$\smash{\Phi_j^{(i)} \in \bbZ[T]}$
when~$j \leq 0$.
\item
For~$X$
projective and
$i$
even, we also have
\begin{equation}
\label{null}
\Phi_j^{(i)} (q^{i/2 - j}) = 0 \, .
\end{equation}
Indeed,~there is the
$\Gal(\overline\bbF_{\!q}/\bbF_{\!q})$-invariant
cycle given by the intersection of
$i/2$~hyperplanes.
The~cycle map~\cite[Cycle, Th\'eor\`eme~2.3.8.iii)]{SGA41/2} yields a non-trivial Galois invariant element
of~$H^i_\et(X_{\overline\bbF_{\!q}}\!, \bbQ_l(i/2))$.
\end{iii}
\end{rems}

\begin{remark}
Consider~the case
that~$i = 1$.
Then
$\smash{N = \dim H^1_\et(X_{\overline\bbF_{\!q}}\!, \bbQ_l)}$
is always even~\cite[Corollaire~(4.1.5)]{De2}. On~the other hand, let a polynomial~$\Phi \in \bbZ[T]$
be given that is of even degree and fulfills assertions a), b), and~c). Then,~by the main theorem of T.\,Honda~\cite{Ho}, there exists an abelian
variety~$A$
such that the eigenvalues
of~$\Frob$
on~$\smash{H^1_\et(A_{\overline\bbF_{\!q}}\!, \bbQ_l)}$
are exactly the zeroes
of~$\Phi$.
One~may enforce that the characteristic polynomial is a power
of~$\Phi$,
and, typically,
$\Phi$~itself
may be~realized.
\end{remark}

\begin{ttt}
We~will show in this note that the same is not true in general
for~$i > 1$.
In~fact, for the characteristic polynomial
of\/~$\Frob$
on the middle cohomology of a variety of even dimension, we will establish a further condition, which is arithmetic in nature and independent of~Theorem~\ref{class}, as well as of formula~(\ref{null}).
\end{ttt}

\begin{theorem}
\label{main}
Let\/~$X$
be a smooth proper variety of even
dimension\/~$d$
over a finite
field\/~$\bbF_{\!q}$
of
characteristic\/~$p$
and\/~$\smash{\Phi = \Phi_{d/2}^{(d)} \in \bbQ[T]}$
be the characteristic polynomial
of\/~$\Frob$
on\/~$\smash{H^d_\et(X_{\overline\bbF_{\!q}}\!, \bbQ_l(d/2))}$.
Put\/~$N := \deg \Phi$.\smallskip

\noindent
Then\/~$(-2)^N \Phi(-1)$
is a square or\/
$p$~times
a square
in\/~$\bbQ$.
\end{theorem}

\begin{remark}
For~$X$
a surface, this result may be deduced from the Tate conjecture, via the Artin-Tate~formula. Cf.~Proposition~\ref{ArT}, below.
\end{remark}

\begin{ttt}
The~correct exponent
of~$p$,
may, at least
for~$p \neq 2$,
be described as~follows.\smallskip

\noindent
{\bf Definition.}
We~put
$$a(X) := \sum_{m=0}^{\frac{d}2-1} (\textstyle\frac{d}2 - m)h'_{d-m,m} \, ,$$
for~$(h'_{d-m,m})_m$
the abstract Hodge numbers
of~$X$
in
degree~$d$\/
\cite[Section\,4]{Maz1}.
\end{ttt}

\begin{rems}
\label{techn}
\begin{abc}
\item
Recall~that the abstract Hodge numbers are defined as~\mbox{follows}. The~crystalline cohomology
groups~$H^i(X/W)$
are finitely generated
\mbox{$W$\!-mod}\-ules,
for
$W := W(\bbF_{\!q})$~the
Witt~ring. They~are acted upon by the absolute Frobenius
$\F$,
the corresponding map is only
\mbox{$\F$-semilinear}~\cite[Expos\'e~I, 2.3.5]{Cham}.

Put~$H := H^i(X/W)/\tors$.
Then,~as
$\F\colon H \to H$
is injective,
$H/\F H$
is a
\mbox{$W$\!-mod}\-ule
of finite~length. By~the classical invariant factor theorem, there is a unique sequence of integers such that
$\smash{H/\F H \cong \!\bigoplus\limits_{m > 0}\! (W/p^mW)^{h'_{i-m,m}}}$.
Finally,~one defines
$\smash{h'_{i,0} := \rk_W \!H - \sum\limits_{m > 0} \!h'_{i-m,m}}$.\smallskip

Observe~that
$a(X)$
is a geometric~quantity. It~depends only on the base
extension~$\smash{X_{\overline\bbF_{\!q}}}$.
\item\looseness-1
Suppose that
$X$
is such that all
$H^i(X/W)$
are torsion-free and that the conjugate spectral sequence
$\smash{E_2^{jm} := H^j(X,\calH^m(\Omega^\bullet_{X/\bbF_{\!q}})) \Longrightarrow H^i_{\rm dR}(X)}$
degenerates
at~$E_2$.
Then
$h'_{i-m,m} = h_{i-m,m} \,(= \dim H^{i-m} (X, \Omega_X^m))$~\cite[Lemma~8.32]{BO}.

For~complete intersections, both assumptions hold~(\cite[Expos\'e~XI, Th\'e\-o\-r\`eme~1.5]{SGA7} together with \cite[Lemma~8.27.2]{BO}). Moreover,~the second assumption is automatically fulfilled when
$\dim X \leq p$
and
$X$
lifts
to~$W$
(\cite[Corollaire~2.4]{DI} and \cite[Lemma 8.27.2]{BO}).
\item
Suppose that
$X$
is of Hodge-Witt type in
degree~$d$,
i.e., that the Serre cohomology groups
$H^j(X,W\!\Omega_X^m)$~\cite{Se}
are finitely generated
\mbox{$W$\!-mod}\-ules
for
$j+m = d$.
The~assertion of Theorem~\ref{main2} below may then be formulated entirely in terms of the characteristic
polynomial~$\Phi$.
In~fact, denote the zeroes of
$\Phi$
by~$z_1, \ldots, z_N$.
Then
$$a(X) = -\!\!\!\!\sum\limits_{\nu_q(z_i) < 0} \!\!\!\!\!\nu_q(z_i)$$
[Corollary~\ref{HodgeWitt}]. This~case includes all varieties that are ordinary in
degree~$d$
\cite[D\'e\-fi\-ni\-tion~IV.4.12]{IR}.
\end{abc}
\end{rems}

\begin{theorem}
\label{main2}
Let\/~$X$
be a smooth proper variety of even
dimension\/~$d$
over a finite
field\/~$\bbF_{\!q}$
of
characteristic\/~$p \neq 2$
and\/~$\smash{\Phi = \Phi_{d/2}^{(d)} \in \bbQ[T]}$
be the characteristic polynomial
of\/~$\Frob$
on\/~$\smash{H^d_\et(X_{\overline\bbF_{\!q}}\!, \bbQ_l(d/2))}$.
Put\/~$N := \deg \Phi$.\smallskip

\noindent
Then\/
$(-2)^N q^{a(X)} \Phi(-1)$
is a square
in\/~$\bbQ$.
\end{theorem}

\begin{rems}
\begin{iii}
\item
It~seems not unlikely that an analogous result is true in
characteristic~$2$,
too. The difficulties that arise may well be of purely technical nature. Cf.~Remark~\ref{rem_char2}, below.
\item
The~assertions may easily be formulated for an arbitrary Tate~twist. But~this does not lead to anything~new. Also,~the twist
by~$d/2$
appears to be natural as the operation
of~$\Frob$
on~$\smash{H^d_\et(X_{\overline\bbF_{\!q}}\!, \bbQ_l(d/2))}$
is orthogonal with respect to a
$\bbQ_l$-valued
symmetric, bilinear~pairing.
\end{iii}
\end{rems}

\begin{rems}
\begin{iii}
\item
One reason for our interest in these conditions is that they simplify the actual computation of the characteristic polynomial
$\Phi$
for a given
variety~$X$.
Such~computations usually involve point counting
on~$X$
over several extensions of the base field, cf.~Examples~\ref{k3} and~\ref{cub4f}. The~conditions given redundantise the most expensive counting~steps.
\item
Applications of the characteristic polynomials
$\Phi$
to varieties in characteristic zero include the determination of the N\'eron-Severi~rank. The~interested reader might consult the articles \cite{vL}, \cite{EJ2}, and~\cite{EJ3}.
\end{iii}
\end{rems}

\noindent
{\bf Acknowledgements.}
We wish to thank Christian Liedtke (M\"unchen), who read an earlier version of this article, for many valuable suggestions.

\section{The proofs}
\label{sec_proo}

\begin{ttt}
{\bf Proof of Theorem~\ref{class}.\;}
a)
This~was first proven by P.~Deligne in~\cite[Th\'e\-o\-r\`eme~(1.6)]{De1} for the projective case and later in \cite[Corollaire~(3.3.9)]{De2}, in~general. The~assertion was formulated by A.\,Weil as a part of his famous conjectures.\smallskip

\noindent
b)
If~$X$
is projective then, by the Hard Lefschetz theorem~\cite[Th\'eor\`eme~(4.1.1)]{De2} and Poincar\'e duality, there is a non-degenerate~pairing
$$\smash{H^i_\et(X_{\overline\bbF_{\!q}}\!, \bbQ_l(j)) \times H^i_\et(X_{\overline\bbF_{\!q}}\!, \bbQ_l(j)) \to \bbQ_l(2j - i)}$$
that is compatible with the operation
of~$\Frob$.
It~is alternating as
$i$
is~odd. The~assertion follows directly from~this. Cf.~the remarks after~\cite[Corollaire~(4.1.5)]{De2}.
The proper non-projective case has only recently been settled by J.~Suh \cite[Corollary~2.2.3]{Su}.\smallskip

\noindent
c)
As~$\Frob$
operates
on~$\smash{H^i_\et(X_{\overline\bbF_{\!q}}\!, \bbZ_l(j))}$,
the
\mbox{$l$-adic}
valuations of the eigenvalues are clearly non-negative. Poincar\'e~duality implies the assertion, cf.~\cite[(2.4)]{De1}.\smallskip

\noindent
d)
This~statement was originally known as Katz's~conjecture. The~usual formulation is that the Newton polygon
of~$\smash{\Phi_0^{(i)}}$
lies above the Hodge polygon of
weight~$i$.
Proofs~are due to B.\,Mazur~\cite{Maz2} and A.\,Ogus~\cite[Theorem~8.39]{BO}.
\eop
\end{ttt}

\begin{nota}
For~$R$
an integral domain,
$K$
its field of fractions, and
$H$
an
\mbox{$R$-module},
we will
write~$H_K := H \!\otimes_R\! K$.
\end{nota}

\begin{lemma}
\label{Cass}
Let\/~$R$
be a principal ideal domain,
$K$
its field of fractions, and\/
$H$~a
free\/
\mbox{$R$-module}
of finite rank, equipped with a perfect,
\mbox{$R$-bilinear},
symmetric pairing\/
$H \times H \to R$.
Denote~its\/
\mbox{$K$-bilinear}
extension by\/
$\langle.\,,.\rangle \colon H_K \times H_K \to K$.\smallskip

\noindent
Furthermore,~let a\/
\mbox{$K$-linear}
map\/
$\sigma\colon H_K \to H_K$
be given that is orthogonal with respect to the~pairing,
i.e.,~$\langle \sigma(x), \sigma(y) \rangle = \langle x, y \rangle$
for
every\/~$x,y \in H_K$.\smallskip

\noindent
Put\/
$$B_0(H) := [H/[H \cap (1-\sigma)H]]_\tors \, . \medskip$$

\begin{abc}
\item
Then~there is a non-degenerate, skew-symmetric\/
\mbox{$R$-bilinear}
pairing
$$(.\,,.) \colon B_0(H) \times B_0(H) \to K/R \, .$$
\item
Suppose\/
$\cha K \neq 2$
and\/
$\langle x, x \rangle \in 2R$
for every\/
$x \in H \cap (1 - \sigma)H_K$.
Then\/~$(.\,,.)$
is~alternating. In~particular, the length
of\/~$B_0(H)$
is~even.
\end{abc}
\end{lemma}

\begin{rems}
\begin{iii}
\item
Observe that
$x \in H$
represents an element
of~$B_0(H)$
if and only
if~$x \in H \cap (1-\sigma)H_K$.
\item
For~$x,y \in H_K$
arbitrary, one has
$\langle (1-\sigma)x, \sigma y \rangle = -\langle x, (1-\sigma)y \rangle$.
In~particular, as
$\sigma\colon H_K \to H_K$
is bijective,
$x \in \ker(1-\sigma)$
if and only if
$x \in ((1-\sigma)H_K)^\perp$.
I.e.,~$(1-\sigma)H_K$
is exactly the set of all elements perpendicular to the
eigenspace~$H_{K,1}$.
This~fact is rather obvious, let us nevertheless emphasize that it is true whether
$\sigma$
is semisimple or~not.
\end{iii}
\end{rems}

\begin{ttt}
{\bf Proof of Lemma~\ref{Cass}.\;}
a)
{\em Definition.}
The~pairing is defined as~follows.
For~$a,b \in B_0(H)$,
choose representatives
$x,y \in H$.
Let~$y' \in H_K$
be such that
$y = (1-\sigma)y'$.
Then~$(a,b) := \langle x, y' \rangle \bmod R$.\smallskip

\noindent
{\em Well-definedness.}
For~two representatives
$x_1, x_2 \in H$,
we have
$x_1 - x_2 = (1-\sigma)v$
for some
$v \in H$.
Thus,
\begin{align*}
\langle x_1 - x_2, y' \rangle = \langle (1-\sigma)v, y' \rangle = \langle (1-\sigma)\sigma v, \sigma y' \rangle & = - \langle \sigma v, (1 - \sigma)y' \rangle \\
 & = - \langle v - x_1 + x_2, y \rangle \in R \, ,
\end{align*}
as both sides are
in~$H$.
On~the other hand, for two representatives
$y'_1, y'_2 \in H_K$,
we have
$(1-\sigma)(y'_1 - y'_2) = (1-\sigma)w \in H$
for a suitable
$w \in H$.
Therefore,
$\langle x, y'_1 - y'_2 \rangle = \langle x, w \rangle + \langle x, k \rangle$
for some
$k \in H_{K,1}$.
The~first summand is
in~$R$,
as both sides are elements
of~$H$.
The~second summand vanishes, since
$x \in (1-\sigma)H_K$.
It~is clear that
$(.\,,.)$
is
\mbox{$R$-bilinear}.\smallskip

\noindent
{\em Non-degeneracy.}
For~$0 \neq b \in B_0(H)$,
one has a
representative~$y$
and some
$y' \in H_K$
such that
$y = (1-\sigma)y'$.
As~$y \not\in (1-\sigma)H$,
we see
$y' \not\in H + H_{K,1}$.
The~goal is to find some
$x \in H \cap (1-\sigma)H_K$
such that
$\langle x, y' \rangle \not\in R$.

For~this, we observe that the perfect pairing induces an isomorphism
$$\smash{H_K \stackrel{\cong}{\longrightarrow} \Hom(H,K)} \, .$$
Under this map,
$H_{K,1} \cong \Hom(H/[H \cap (1-\sigma)H_K], K)$.
Furthermore,
$$H + H_{K,1} \cong \{ \alpha \in \Hom(H,K) \mid \alpha(H \cap (1-\sigma)H_K) \subseteq R \} \, .$$
Indeed,~as
$H \cong \Hom(H,R)$,
the inclusion
``$\subseteq$''
is~obvious. The~other inclusion follows from the fact that
$H \cap (1-\sigma)H_K$
is a direct summand
of~$H$.
The~homomorphism
$\alpha|_{H \cap (1-\sigma)H_K} \colon H \cap (1-\sigma)H_K \to R$
may thus be extended to a homomorphism
$\alpha' \colon H \to R$,
corresponding to an element
of~$H$.
The~difference
$\alpha - \alpha'$
vanishes
on~$H \cap (1-\sigma)H_K$,
hence is defined by an element
of~$H_{K,1}$.

Now,~as
$y' \not\in H + H_{K,1}$,
the corresponding homomorphism does not send
$H \cap (1-\sigma)H_K$
to~$R$.
I.e.,~there is some
$x \in H \cap (1-\sigma)H_K$,
not mapped
to~$R$.
This~is exactly our~claim.\smallskip

\noindent
{\em Skew-symmetry.}
Let
$a,b \in B_0(H)$
and choose representatives
$x,y \in H$.
There~are
$x',y' \in H_K$
such that
$x = (1-\sigma)x'$
and
$y = (1-\sigma)y'$.
Then
$$\langle x, y' \rangle + \langle x', y \rangle = \langle x' - \sigma x', y' \rangle + \langle x', y' - \sigma y' \rangle = \langle x' - \sigma x', y' - \sigma y' \rangle = \langle x, y \rangle \in R \, .$$

\noindent
b)
For~$a \in B_0(H)$
choose a representative
$x \in H$
and
$x' \in H_K$
such that
$x = (1-\sigma)x'$.
Then
\begin{align*}
2 (a,a) = 2\langle x, x' \rangle = [\langle x', x' \rangle - \langle \sigma x', x' \rangle] + {} & [\langle \sigma x', \sigma x' \rangle - \langle x', \sigma x' \rangle] \\
& = \langle x' - \sigma x', x' - \sigma x' \rangle = \langle x, x \rangle \in 2R \, ,
\end{align*}
hence
$(a, a) = 0$.\smallskip

\noindent
{\em Evenness of~$\ell(B_0(H))$.}
First~note that
$B_0(H)$
is a finitely generated torsion module over a principal ideal~domain. This~implies that
$\ell(B_0(H))$
is~finite.

It~is known~\cite[Exercise~10.20]{Si} that, for such a
module~$M$,
the existence of a non-degenerate alternating pairing implies that
$\ell(M)$
is~even. The~argument is essentially as~follows.

Assume~the contrary and let
$M$
be a counterexample of minimal~length.
Choose~$x \in M$
such that
$Rx \subseteq M$
is of length~one.
I.e.,~$\frakm := \Ann_R(x)$
is a maximal~ideal. Then
$L_x := \langle x,.\, \rangle \colon M \to K/R$
is a nonzero linear form, the image of which
is~$\frakm^{-1}/R \cong R/\frakm$.
Clearly,~one
has~$x \in \ker L_x$.

There~is an alternating form induced
on~$[\ker L_x]/(x)$,
which is again non-de\-gen\-er\-ate.
Moreover,~$\ell(\ker L_x) = \ell(M)-1$
and
$\ell([\ker L_x]/(x)) = \ell(M)-2$.
This~is a contradiction.
\eop
\end{ttt}

\begin{rems}
\begin{iii}
\item
If~$2$~is
a unit
in~$R$
then the assumption of~b) is automatically~fulfilled.
\item
When
$R = \bbZ_l$,
$\sigma(H) \subseteq H$,
and
$\sigma$
is semisimple
at~$1$,
this~result was proven by Yu.\,G.~Zarhin in~\cite[3.3 and~Lemma~3.4.1]{Zar}. It~is implicitly contained in the work of J.\,W.\,S.~Cassels~\cite{Ca1}.
\item
Most~of our applications will be based on the following~corollary.
\end{iii}
\end{rems}

\begin{corollary}
\label{det}
Let\/~$(R,\nu)$
be a normalized discrete valuation ring,
$k$~its
residue field, which we assume to be finite of characteristic\/
$\neq\! 2$,
$K$
the field of fractions, and\/
$H$~a
free\/
\mbox{$R$-module}
of finite~rank. Suppose~there is a perfect,
\mbox{$R$-bilinear},
symmetric pairing\/
$H \times H \to R$
and denote its
\mbox{$K$-bilinear}
extension by\/
$\langle.\,,.\rangle \colon H_K \times H_K \to K$.\smallskip

\noindent
Furthermore,~let a\/
\mbox{$K$-linear}
map\/
$\sigma\colon H_K \to H_K$
be given that is orthogonal with respect to the pairing and such that\/
$1$
is not among its~eigenvalues.

\begin{abc}
\item
Then\/
$\ell(\sigma(H)+H/H) + \nu(\det(1-\sigma))$
is~even.
\item
In~particular, if\/
$\sigma(H) \subseteq H$
then\/
$\nu(\det(1-\sigma))$
is~even.
\end{abc}\smallskip

\noindent
{\bf Proof.}
{\em
a)
We~write
$q := \#k$.
According~to the definition, the modulus \cite[section~14.3]{Di} of the map
$(1-\sigma)\colon H_K \to H_K$
may be computed as
$q^{\ell(M/H) - \ell(M/(1-\sigma)H)}$
for every
\mbox{$R$-module}
$M \supseteq H, (1-\sigma)H$
that is chosen such that the two lengths are~finite. Moreover,~by \cite[Exercise~14.3.6]{Di}, that modulus is equal to
$q^{-\nu(\det(1-\sigma))}$.
Consequently,
$$\nu(\det(1-\sigma)) = \ell([\sigma(H) + H]/(1-\sigma)H) - \ell(\sigma(H) + H/H) \, .$$

On~the other hand, all the assumptions of Lemma~\ref{Cass}.b) are~fulfilled.
As~$1$
is not an eigenvalue
of~$\sigma$,
$H/[H \cap (1-\sigma)H]$
is purely~torsion. Thus,~we have that
$\ell(H/[H \cap (1-\sigma)H])$
is~even. Furthermore,
$$H/[H \cap (1-\sigma)H] \cong [H + (1-\sigma)H] / (1-\sigma)H = [\sigma(H) + H] / (1-\sigma)H \, .$$
I.e.,~$\ell([\sigma(H) + H] / (1-\sigma)H)$
is even,~too. The~assertion~follows.\smallskip

\noindent
b)~is an immediate consequence of~a).
}
\eop
\end{corollary}

\begin{ttt}
In~order to illustrate the strength of Corollary~\ref{det}, let us show an application to modules of rank two, the smallest non-trivial~case. The~fact obtained belongs to the not-so-well-known results on real quadratic number~fields. Cf.~\cite[p.\,118]{Zag}.\smallskip

\noindent
{\bf Corollary.}
{\em
Let\/~$K = \bbQ(\sqrt{d})$
be a quadratic number field
and\/~$\varepsilon \in K$
a unit of
norm\/~$(+1)$.
Then\/~$\N(1 - \varepsilon) = 2 - \tr(\varepsilon)$
is a product of some primes dividing the discriminant, a perfect square, and, possibly, a factor\/~$2$
and a minus~sign.\smallskip

\noindent
{\bf Proof.}
{\em
As~$K$
is a quadratic number field, the norm
$\N\colon \calO_K \to \bbZ$
is a quadratic~form. The~multiplication map
$\cdot\varepsilon\colon \calO_K \to \calO_K$
is compatible with this form and, therefore, orthogonal with respect to the symmetric, bilinear form
$\langle.\,,.\rangle \colon \calO_K \times \calO_K \to \bbZ$
associated
to~$\N$.

The~same is true for the corresponding
\mbox{$\bbZ_l$-valued}
pairings between the
\mbox{$l$-adic}
completions
of~$\calO_K$,
for
$l$
any prime~number. As~these pairings are perfect as long as
$l$
does not divide the discriminant
of~$K$,
the assertion follows from Corollary~\ref{det}.b).
}}
\eop
\end{ttt}

\begin{ttt}
{\bf Proof of Theorem~\ref{main}.\;}
We~clearly have that
$(-2)^N \Phi(-1) \in \bbQ$.
Furthermore,~we may assume
that~$(-1)$
is not among the zeroes
of~$\Phi$
as, otherwise, the assertion is true,~trivially.

Then~$(-2)^N \Phi(-1) > 0$.
Indeed,~as
$\Phi \in \bbR[T]$
and there is no real zero different
from~$1$,
$$(-2)^N \Phi(-1) = 2^N (1 + z_1) \cdot\ldots\cdot (1 + z_N)$$
is the product of several factors of the form
$z\overline{z} = |z|^2$
for~$z \in \bbC$
and some factors that are equal
to~$2$.
To~prove the assertion, we will show that
$(-2)^N \Phi(-1)$
is of even
\mbox{$l$-adic}
valuation for every prime
number~$l \neq p$.\smallskip

\noindent
Put~$\smash{H := H^d_\et(X_{\overline\bbF_{\!q}}\!, \bbZ_l(d/2)) / \tors}$.
By~Poincar\'e duality \cite[Exp.~XVIII, formule (3.2.6.2)]{SGA4}, cf.~\cite[Chap.~6, Sec.~2, Theorem~18 and Chap.~5, Sec.~5, Theorem~3]{Sp}, the bilinear~pairing
$$\langle.\,,.\rangle \colon H \times H \longrightarrow H^{2d}_\et(X_{\overline\bbF_{\!q}}\!, \bbZ_l(d)) \stackrel{\cong}{\longrightarrow} \bbZ_l \, ,$$
given by cup~product and trace~map, is~perfect.
As~$d$
is even, it is symmetric,~too. The~operation
of~$\Frob$
on~$H_{\bbQ_l}$
is orthogonal with respect to the
\mbox{$\bbQ_l$-linear}
extension of this~pairing.\medskip

\noindent
{\em First case.}
$l \neq 2$.\smallskip

\noindent
The~operation
of~$(-\Frob)$
is orthogonal with respect to the pairing,~too.
As~$1$
is not among its eigenvalues, Corollary~\ref{det}.b) shows that
$\nu_l(\det (1+\Frob)) = \nu_l(\Phi(-1))$
is~even.\medskip

\noindent
{\em Second case.}
$l = 2$.\smallskip

\noindent\looseness-1
Here,~the argument is a bit more~involved. First,~note that the tangent sheaf
$\calT_X$
of~$X$
is defined over the base
field~$\bbF_{\!q}$.
This~shows that the Chern classes
$\smash{c_i(\calT_X) \in H^{2i}_\et(X_{\overline\bbF_{\!q}}\!, \bbZ_2(i))}$
are invariant
under~$\Frob$.
We~therefore see from Lemma \ref{Wu} that there is a
$\Frob$-invariant
element
$\smash{\omega \in H^d_\et(X_{\overline\bbF_{\!q}}\!, \bbZ_2(d/2))}$
such that
$\langle \omega, x \rangle + \langle x,x \rangle \in 2\bbZ_2$
for every
$\smash{x \in H^d_\et(X_{\overline\bbF_{\!q}}\!, \bbZ_2(d/2))}$.

For~$x \in (1 - \Frob)H_{\bbQ_2}$,
the fact that
$\omega$
is
\mbox{$\Frob$-invariant}
implies~$\langle \omega, x \rangle = 0$.
Hence,
$\langle x, x \rangle \in 2\bbZ_2$
for~$x \in H \cap (1 - \Frob)H_{\bbQ_2}$.
According~to Lemma~\ref{Cass}.b), 
$[H/(1-\Frob)H]_\tors$
is of even~length.

An~application
to~$\smash{X_{\bbF_{\!q^2}}}$
shows that
$[H/(1-\Frob^2)H]_\tors$
is of even length,~too. Lemma~\ref{Hilf} now yields the~assertion.
\eop
\end{ttt}

\begin{lemma}
\label{Wu}
For~every even\/
$d > 0$,
there exists a polynomial\/
$P_d \in \bbZ[T_0,\ldots, T_{d/2}]$,
weighted homogeneous of degree\/
$d$
for\/
$T_i$
of
weight\/~$2i$,
such that the following is~true.\smallskip

\noindent
For~a smooth projective
variety\/~$X$
of
dimension\/~$d$
over a finite
field\/~$\bbF_{\!q}$
of
characteristic\/~$\neq 2$,
put\/
$\smash{\omega := P_d(c_0(\calT_X), \ldots c_{d/2}(\calT_X)) \in H^d_\et(X_{\overline\bbF_{\!q}}\!, \bbZ_2(d/2))}$
for\/
$\calT_X$
the tangent sheaf and\/
$c_i$
the\/
$i$-th
Chern~class. Then
$$\langle \omega, x \rangle + \langle x,x \rangle \in 2\bbZ_2$$
for
every\/~$\smash{x \in H^d_\et(X_{\overline\bbF_{\!q}}\!, \bbZ_2(d/2))}$.\smallskip

\noindent
{\bf Proof.}
{\em
This~is a standard result in topology, cf.~\cite[\S\S8 and~11]{MS}, and carries over to the \'etale site. Denote~by
$\smash{\varepsilon\colon H^{2*}_\et(X_{\overline\bbF_{\!q}}\!, \bbZ_2(*)) \longrightarrow H^{2*}_\et(X_{\overline\bbF_{\!q}}\!, \bbZ/2\bbZ)}$
the reduction~map. Furthermore,~for simplicity, write
$c_i := c_i(\calT_X) \in H^{2i}_\et(X_{\overline\bbF_{\!q}}\!, \bbZ_2(i))$
for the Chern~classes.

For~$k \in \bbN_0$,
let
$\smash{\nu_{2k} \in H^{2k}_\et(X_{\overline\bbF_{\!q}}\!, \bbZ/2\bbZ)}$
be the
\mbox{$2k$-th}
Wu class
of~$X$~\cite[p.\,578]{Ur}.
If~$X$
is of dimension
$d$
then, according to the very definition of the Wu~class,
$$\nu_d \cup x + x \cup x = \Sq^d(x) + \Sq^d(x) = 0$$
for every
$\smash{x \in H^d_\et(X_{\overline\bbF_{\!q}}\!, \bbZ/2\bbZ)}$~\cite[Prop.~2.2.(2)]{Ur}.
We~will inductively construct polynomials
$P_{2k} \in \bbZ[T_0,\ldots, T_k]$
such that
$\varepsilon(P_{2k}(c_0(\calT_X), \ldots c_{d/2}(\calT_X))) = \nu_{2k}$.

For~every
$k \in \bbN_0$,
there is the formula of Wu~\cite[Proposition~0.5]{Ur},
$$\Sq^0(\nu_{2k}) + \Sq^2(\nu_{2k-2}) + \cdots + \Sq^{2k}(\nu_0) = \varepsilon(c_k) \, .$$
Moreover,~for the Steenrod squares of the Chern classes, there are the formulas
$$\Sq^{2j}(\varepsilon(c_i)) \!=\! \varepsilon(c_j)\varepsilon(c_i) + \left(\! \atop{2j-2i}2 \!\right)\! \varepsilon(c_{j-1})\varepsilon(c_{i+1}) + \cdots + \left(\! \atop{2j-2i}{2j} \!\right)\! \varepsilon(c_0)\varepsilon(c_{i+j}) .$$
Indeed,~these follow in a purely formal manner from the definitions of Chern classes and Steenrod squares, cf.~\cite[Problem~8-A]{MS}.

As~$\Sq^0 = \id$,
Wu's formula implies that we may choose
$P_0(T_0) := T_0$.
Furthermore,~having
$P_0, \ldots, P_{2k-2}$
already constructed, it shows~that
\begin{align*}
\nu_{2k} = \varepsilon(c_k) - \sum_{i=0}^k \Sq^{2i}(\nu_{2k-2i}) &= \varepsilon(c_k) - \sum_{i=0}^k \Sq^{2i}(\varepsilon(P_{2k-2i}(c_0, \ldots, c_{k-i}))) \\
&= \varepsilon(c_k) - \sum_{i=0}^k P_{2k-2i}(\Sq^{2i}(\varepsilon(c_0), \ldots, \Sq^{2i}(\varepsilon(c_{k-i})))) \, .
\end{align*}
Plugging~into this the formula for the Steenrod squares of the Chern classes, we see that
$\nu_{2k}$
is the reduction of a polynomial expression
in~$c_0, c_1, \ldots, c_k$.
}
\eop
\end{lemma}

\begin{lemma}
\label{Hilf}
Let\/~$(R,\nu)$
be a normalized discrete valuation ring of
characteristic\/~$0$,
$k$~its
residue field, which we assume to be finite,
$K$
the field of fractions, and\/
$H$
a free\/
\mbox{$R$-module}
of finite rank, equipped with a non-degenerate, symmetric\/
\mbox{$K$-bilinear}
pairing\/
$\langle.\,,.\rangle \colon H_K \times H_K \to K$.\smallskip

\noindent
Moreover,~let an\/
\mbox{$R$-linear}
map\/
$\sigma\colon H \to H$
be given that is orthogonal with respect to the pairing and does not have the
eigenvalue\/~$(-1)$.
Suppose~that
$[H/(1 - \sigma)H]_\tors$
and
$[H/(1 - \sigma^2)H]_\tors$
are of even~lengths.\smallskip

\noindent
Then
$$\nu(\det(1+\sigma)) \equiv \nu(2) \!\cdot\! \rk H \pmod 2 \, .$$

\noindent
{\bf Proof.}
{\em
We~will prove this technical lemma in several~steps.\smallskip

\noindent
{\em First step.}
$(1 - \sigma)H/(1 - \sigma^2)H$
is of even~length.

\noindent
Since~$(-1)$
is not an eigenvalue
of~$\sigma$,
the map
$(1+\sigma)\colon H_K \to H_K$
is a~bijection. In~particular,
$(1 - \sigma)H_K = (1 - \sigma^2)H_K$.
Furthermore,~there is the commutative diagram of short exact~sequences
$$
\definemorphism{gleich}\Solid\notip\notip
\xymatrixcolsep{2.7mm}
\xymatrix{
0 \rto & (1 \!-\! \sigma)H/(1 \!-\! \sigma^2)H \rto\ar@{->>}[d] & H/(1 \!-\! \sigma^2)H \rto\ar@{->>}[d] & H/(1 \!-\! \sigma)H \rto\ar@{->>}[d] & 0 \phantom{\, .} \\
0 \rto & 0 \rto & H/H \!\cap\! (1 \!-\! \sigma^2)H_K \rgleich & H/H \!\cap\! (1 \!-\! \sigma)H_K \rto & 0 \, .
}
$$
As~the vertical arrows are surjective, the
\mbox{$9$-lemma}
yields exactness~of
$$0 \longrightarrow (1 - \sigma)H / (1 - \sigma^2)H \longrightarrow [H/(1 - \sigma^2)H]_\tors \longrightarrow [H/(1 - \sigma)H]_\tors \longrightarrow 0 \, .$$
Since~$\ell([H/(1 - \sigma)H]_\tors)$
and
$\ell([H/(1 - \sigma^2)H]_\tors)$
are even, we see that the length
of
$(1 - \sigma)H/(1 - \sigma^2)H$
is even,~too. This~is the~claim.\smallskip\pagebreak[3]

\noindent
{\em Second step.}
$\nu(\det (1+\sigma)|_{(1 - \sigma)H_K})$
is~even.

\noindent
Both,
$(1 - \sigma)H$
and
$(1 - \sigma^2)H$,
are
\mbox{$R$-submodules}
of maximal rank in the
\mbox{$K$-vector}
space
$(1 - \sigma)H_K$.
Moreover,~$(1 - \sigma^2)H = (1 + \sigma)|_{(1 - \sigma)H_K} [(1 - \sigma)H]$.
Hence,~the modulus of
$(1 + \sigma)|_{(1 - \sigma)H_K}$
is
$\smash{q^{-\ell((1 - \sigma)H/(1 - \sigma^2)H)}}$
for
$q := \#k$.

On~the other hand, by~\cite[Exercise~14.3.6]{Di}, this modulus is
$\smash{q^{-\nu(\det (1+\sigma)|_{(1 - \sigma)H_K})}}$.
The~result established in the first step implies the~claim.\smallskip

\noindent
{\em Third step.}
$\nu(\det(1+\sigma)) \equiv [\dim \ker (1 - \sigma)] \nu(2) \pmod 2$.

\noindent
According~to C.~Jordan, we have
$H_K = \ker (1 - \sigma)^r + (1 - \sigma)H_K$
for~$r$~large.
On~$\ker (1 - \sigma)^r$,
the
homomorphism~$\sigma$
has only the
eigenvalue~$1$.
In~particular,
$(1 - \sigma)H_K \subseteq H_K$
has a
complement~$V$
such that
$2$~is
the only eigenvalue
of~$(1+\sigma) |_V$.
Hence,
\begin{align*}
\nu(\det(1+\sigma)) &   =    \smash{\nu(\det(1+\sigma))|_V + \nu(\det(1+\sigma))|_{(1 - \sigma)H_K}} \\
                    & \equiv \nu(\det(1+\sigma))|_V  \hspace{2cm} \pmod 2 \\
                    &   =    \nu(2) [\dim V] \\
                    &   =    \nu(2) [\rk H - \dim (1 - \sigma)H_K] \\
                    &   =    \nu(2) [\dim \ker (1 - \sigma)] \, .
\end{align*}
{\em Fourth step.}
Orthogonality.

\noindent
Finally,~$\sigma$
is orthogonal with respect to a non-degenerate bilinear~form. In~this~situation, it is well-known \cite[Example~2.6.C.i.b)]{Wa} (see also \cite[\S1.A1, Lemma~4.ii)]{St})~that
$$\dim \ker (1 - \sigma) \equiv \dim \ker (1 - \sigma)^r \pmod 2$$
for~$r \gg 0$.

The~right hand side is the algebraic multiplicity of the
eigenvalue~$1$.
As~the eigenvalues of an orthogonal matrix come in pairs
$(z, \frac1z)$
and, by our assumption,
$(-1)$
is not an eigenvalue, it has the same parity
as~$\rk H$.
The~assertion~follows.
}
\eop
\end{lemma}

\begin{rems}
\begin{iii}
\item
There~is a conjecture of J.-P.~Serre that the operation
of~$\Frob$
on
\mbox{$l$-adic}
cohomology is always semisimple.
Then,~as the eigenvalues come in pairs
$(z, \frac1z)$,
the congruence
$\dim \ker (1 - \Frob) \equiv N \pmod 2$
is~clear. Thus,~conjecturally, the argument on the sizes of the Jordan blocks is not necessary in the application to
\mbox{$l$-adic}~cohomology.
\item
Suppose~that
$q = p^k$
for a prime
$p \neq 2$
and let
$X$
be a surface such that the canonical sheaf
$K \in \Pic(X_{\overline\bbF_{\!q}})$
is divisible
by~$2$.
Then~the case
$l=2$
of Theorem~\ref{main} may be treated~directly.

Indeed,~in this situation, Wu's~formula~\cite[Proposition~2.1]{Ur} implies that
$\langle x, x \rangle \in 2\bbZ_2$
for
every~$x \in H$.
Therefore,~by Lemma~\ref{Cass}.b),
$H/(1 + \Frob)H$
is of even~length. Moreover,~the assumption enforces that
$K^2$
is~even.
Hence,
$N = \dim H^2_\et(X_{\overline\bbF_{\!q}}\!, \bbQ_2(1))$
is even, too, by Noether's formula~\cite[I.14]{Bea}.
\end{iii}
\end{rems}

\begin{ttt}
{\bf Proof of Theorem~\ref{main2}.\;}
Here,~according to our assumption, we
have~$p \neq 2$.
Again,~we may assume without restriction
that~$\Phi(-1) \neq 0$.
In~view of Theorem~\ref{main}, it will suffice to prove that
$q^{a(X)} \Phi(-1)$
is of even
\mbox{$p$-adic}~valuation.
Writing~$q = p^k$,
this means that
$ka(X) + \nu_p((1+z_1)\cdots(1+z_N))$
is~even.

For~this, let
$W := W(\bbF_{\!q})$
be the Witt ring and
$K$
its field of~fractions. The~crystalline cohomology groups
$H^d(X/W)$
are finitely generated
\mbox{$W$\!-mod}\-ules,
acted upon semilinearly by the absolute
Frobenius~$\F$.
The~operation
of~$\F^k$
is
\mbox{$W$-linear}
and such that its characteristic polynomial coincides
with~$\smash{\Phi_0^{(d)}}$~\cite{KM}.

Again,~$H := H^d(X/W)/\tors$
is equipped with a natural perfect pairing~\cite[Ch.~VII, Th\'eor\`eme~2.1.3]{Be}
$$\langle.\,,.\rangle \colon H \times H \longrightarrow W \, .$$
The~Frobenius operation is, however, compatible with this pairing only in the sense that
$\langle \F(x), \F(y) \rangle = p^d \langle x, y \rangle$
for
every~$x,y \in H$~\cite[Expos\'e~II, Exemple~1.1.ii)]{Cham}.
The~map
$\sigma := \F\!/p^{d/2} \colon H_K \to H_K$
respects the
pairing~$\langle.\,,.\rangle$.
$H$~is,
in fact, a Dieu\-donn\'e module and thus carries a rich structure~\cite{Man}. We~shall use only a small part of~it. Let~us distinguish two~cases.\smallskip

\noindent
{\em First case.}
$k$ is odd.

\noindent
Observe~that
$H$
is a free module as well over the Witt ring
$W(\bbF_{\!p})$
and
$\sigma$
is
$W(\bbF_{\!p})$-linear.
Clearly,~$\rk_{W(\bbF_{\!p})} H = k \!\cdot\! \rk_W \!H$.
The~eigenvalues
of~$\sigma$,
as a
$W(\bbF_{\!p})$-linear
map, are all the
\mbox{$k$-th}
roots of the zeroes
$z_1, \ldots, z_N$
of~$\Phi$.

Furthermore,~$\sigma$
is orthogonal with respect to the
$\Q_p := \Q(W(\bbF_{\!p}))$-bilinear
extension
$H_K \times H_K \to \Q_p$
of the perfect pairing
$\smash{\tr_{W(\bbF_{\!p^k})/W(\bbF_{\!p})} \circ\, \langle.\,,.\rangle \colon H \times H \to W(\bbF_{\!p})}$.
The~operation of
$(-\sigma)$
is orthogonal with respect to the pairing,~too. Thus,~Corollary~\ref{det} shows that
$\ell_{W(\bbF_{\!p})}(\sigma(H)+H/H) + \nu_{W(\bbF_{\!p})}(\det(1+\sigma))$
is~even. By~Lemma~\ref{ee}, the first summand is equal
to~$ka(X)$.
The~second summand is the
\mbox{$p$-adic}
valuation~of
$$\prod_{i=1}^N \prod_{r^k = z_i} \!(1+r) = \prod_{i=1}^N (1+z_i) \, .$$

\noindent
{\em Second case.}
$k$ is even.

\noindent
Here,~the argument is slightly more~involved. We~first observe that
$H$
is a~free \mbox{module} over the Witt ring
$\smash{W(\bbF_{\!p^{k/2}})}$,
too, and that
$\sigma^{k/2}$
is
$\smash{W(\bbF_{\!p^{k/2}})}$-linear.
Clearly,
$\smash{\rk_{W(\bbF_{\!p^{k/2}})} H  = 2 \!\cdot\! \rk_W \!H}$.
Moreover,~$\sigma^{k/2}$
is orthogonal with respect to the
$\Q_{p^{k/2}} := \Q(W(\bbF_{\!p^{k/2}}))$-bilinear
extension
$H_K \times H_K \to \Q_{p^{k/2}}$
of the perfect pairing
$\smash{\tr_{W(\bbF_{\!p^k})/W(\bbF_{\!p^{k/2}})} \circ\, \langle.\,,.\rangle \colon H \times H \to W(\bbF_{\!p^{k/2}})}$.

Now~choose a unit
$\smash{u \in W(\bbF_{\!p^k})}$
such that
$\smash{\tr_{W(\bbF_{\!p^k})/W(\bbF_{\!p^{k/2}})} (u) = 0}$,
i.e.,~such that its conjugate
is~$(-u)$.
Then
$\smash{\sigma^{k/2}(ux) = -u \sigma^{k/2}(x)}$
for
all~$x \in H$.

We~define a
$\smash{W(\bbF_{\!p^{k/2}})}$-linear
map
$T\colon H_K \to H_K$
by
$$T(x) := \sigma^{k/2} (ux) \, .$$
This~yields
$T \!\circ\! T = -u^2\sigma^k$.
Furthermore,
$T$~is
$\smash{W(\bbF_{\!p^k})}$-semilinear
and one has
$\langle Tx, Ty \rangle = \langle ux, uy \rangle = u^2 \langle x, y \rangle$
for all
$x,y \in H$.

We~see that the eigenvalues
of~$T$,
as a
$\smash{W(\bbF_{\!p^{k/2}})}$-linear
map, are all the square roots of the numbers
$(-u^2 z_1), \ldots, (-u^2 z_N)$,
for
$z_1, \ldots, z_N$
the zeroes
of~$\Phi$.
Moreover,~$[T(H) + H]/H$
is of even
\mbox{$\smash{W(\bbF_{\!p^{k/2}})}$-length}
as it is, in fact, a
\mbox{$\smash{W(\bbF_{\!p^k})}$-module}.

Finally,~put
$\smash{\overline{H} := H \!\otimes_{W(\bbF_{\!p^{k/2}})}\! W(\overline\bbF_{\!p})}$
and extend
$T$
to a
\mbox{$W(\overline\bbF_{\!p})$-linear}
map
$\overline{T}\colon \overline{H} \to \overline{H}$.
Then~$\frac1u \overline{T}$
is~orthogonal. On~the other~hand,
$$\textstyle [\frac1u \overline{T}(\overline{H}) + \overline{H}]/\overline{H} = [\overline{T}(\overline{H}) + \overline{H}]/\overline{H} \, ,$$
as
$\smash{\frac1u}$
is a~unit, and the latter
\mbox{$W(\overline\bbF_{\!p})$-module}
is of even~length. Corollary~\ref{det} shows that
$\nu_{W(\overline\bbF_{\!p})}(\det(1 - \frac1u \overline{T}))$
is~even. This~number is nothing but the
\mbox{$p$-adic}
valuation~of
\begin{align*}
\prod_{i=1}^N \prod_{r^2 = -z_i} \!\!\!(1-r) = \prod_{i=1}^N (1+z_i) \, . \eop
\end{align*}
\end{ttt}

\begin{remark}
\label{rem_char2}
One~might want to prove Theorem~\ref{main2} in
characteristic~$2$
along the same lines as Theorem~\ref{main}. For~this, one would need, at least, the theory of Steenrod squares and Wu classes, as well as Wu's~formula, for crystalline cohomology of varieties in
characteristic~$2$.
It~seems, however, that such a theory is not yet available in the~literature.
\end{remark}

\begin{lemma}
\label{ee}
Let\/~$X$
be a smooth proper variety of even
dimension\/~$d$
over\/~$\bbF_{\!q}$,
$H := H^d(X/W)/\tors$,
and\/
$\sigma = \F\!/p^{d/2}$.
Then\/~$a(X) = \ell_W (\sigma(H)+H/H)$.\smallskip

\noindent
{\bf Proof.}
{\em
First,~observe that
$\sigma(H)$,
being the image of an
\mbox{$\F$-semilinear}
map, is indeed
a~\mbox{$W$\!-mod}\-ule.
Furthermore,~we have
$\smash{H/\F H \cong \bigoplus_{m > 0} (W/p^mW)^{h'_{d-m,m}}}$.
Hence,~there is a basis
of~$H$
such that, under the corresponding isomorphism
$H \cong W^N$,
one has
$\smash{\F H \cong \bigoplus_{m \geq 0} p^mW^{h'_{d-m,m}}}$.
Therefore,~$\smash{\sigma(H) \cong \!\bigoplus\limits_{m \geq 0}\!\! p^{m-d/2} W^{h'_{d-m,m}}}$.
The~assertion~follows.
}
\eop
\end{lemma}

\begin{corollary}
\label{HodgeWitt}
Let\/~$X$
be a smooth proper variety of even
dimension\/~$d$
over\/~$\bbF_{\!q}$.
Suppose that\/
$X$
is of Hodge-Witt type in
degree~$d$,
i.e., that the Serre cohomology groups\/
$H^j(X,W\!\Omega_X^m)$
are finitely generated\/
\mbox{$W$\!-mod}\-ules
for~$j+m = d$.
Then
$$a(X) = -\!\!\!\!\sum\limits_{\nu_q(z_i) < 0} \!\!\!\!\!\nu_q(z_i) \, .$$
{\bf Proof.}
{\em
We~have
$\smash{H^d(X/W) \cong \bigoplus_m H^{d-m,m}}$
for
$H^{d-m,m} := H^{d-m}(X,W\!\Omega_X^m)$,
as is shown
in~\cite[Th\'eor\`eme~IV.4.5]{IR}.
On~$H^{d-m,m}$,
$\F$
operates as
$p^m F$
for~$F$
the usual Frobenius on Serre~cohomology.
Thus,~$\sigma$
acts
as~$\smash{p^{m - d/2}F}$.
For~$m \geq d/2$,
this ensures that the corresponding summand is mapped
to~$H$.

Thus,~assume
that~$m < d/2$.
On~Serre cohomology, there is a second operator, the
Verschiebung~$V$,
such that
$FV = p$.
Hence~$\sigma p^{d/2-m-1} V = \id$,
implying
$$H^{d-m,m} \!\otimes_W\! \Q(W) \supseteq \sigma(H^{d-m,m}) \supseteq H^{d-m,m} \, .$$
Lemma~\ref{ee} shows that
$\smash{a(X) = - \nu_q(\det(\sigma|_{\!\!\!\!\bigoplus\limits_{m < d/2} \!\!\!\!\!H^{d-m,m}}\!))}$,
which is equivalent to the~assertion.
}
\eop
\end{corollary}

\section{An application and examples}

\subsection{An application to the odd-dimensional supersingular~case}

\begin{lemma}
\label{comb}
Let\/~$p$
be a prime and\/
$d_1$
and\/~$d_2$
two odd~integers. Moreover,~let\/
$\Phi_1 \in \bbQ[T]$
and\/
$\Phi_2 \in \bbQ[T]$
be polynomials of even
degrees\/~$N_1$
and\/~$N_2$
that fulfill the functional~equations
$$\Phi_1(p^{d_1}/T) = \frac{p^{d_1 N_1 / 2}}{T^{N_1}} \Phi_1(T)
\quad{\text and~}\quad
\Phi_2(p^{d_2}/T) = \frac{p^{d_2 N_2 / 2}}{T^{N_2}} \Phi_2(T) \, .$$
For\/~$z_i^{(1)}$
the zeroes
of\/~$\Phi_1$
and\/
$z_j^{(2)}$
the zeroes
of\/~$\Phi_2$,
let\/
$\Phi$
be the monic \mbox{polynomial} with the
zeroes\/~$\smash{z_i^{(1)} \!z_j^{(2)} \!/ p^{\frac{d_1 + d_2}2}}$.\smallskip

\noindent
Then\/
$p^{N_1N_2/4} \Phi(-1)$
is a square
in\/~$\bbQ$.\medskip

\noindent
{\bf Proof.}
{\em
The~assumption implies that the zeroes come in pairs with
products~$p^{d_1}$
and~$p^{d_2}$,~respectively.
For~two pairs of zeroes, the corresponding four zeroes
of~$\Phi$~are
$$\frac{z_i^{(1)}}{p^{d_1 / 2}} \!\cdot\! \frac{z_j^{(2)}}{p^{d_2 / 2}},\quad
\frac{z_i^{(1)}}{p^{d_1 / 2}} \!\cdot\! \frac{p^{d_2 / 2}}{z_j^{(2)}},\quad
\frac{p^{d_1 / 2}}{z_i^{(1)}} \!\cdot\! \frac{z_j^{(2)}}{p^{d_2 / 2}},\quad
\textrm{\upshape and}~ \quad
\frac{p^{d_1 / 2}}{z_i^{(1)}} \!\cdot\! \frac{p^{d_2 / 2}}{z_j^{(2)}} \, .$$
We~now observe the identity
\begin{equation}
\label{quadm}
\textstyle (-1 - \frac1p u_1u_2)(-1 - u_1/u_2)(-1 - u_2/u_1)(-1 - p/u_1u_2) = p(\frac{u_1}p + \frac{u_2}p + 1/u_1 + 1/u_2)^2 \, ,
\end{equation}
which applies, since the four zeroes may rationally be written as
$$\frac1p \frac{z_i^{(1)}}{p^{(d_1-1) / 2}} \!\cdot\! \frac{z_j^{(2)}}{p^{(d_2-1) / 2}}, \,
\frac{z_i^{(1)}}{p^{(d_1-1) / 2}} \!\cdot\! \frac{p^{(d_2-1) / 2}}{z_j^{(2)}}, \,
\frac{p^{(d_1-1) / 2}}{z_i^{(1)}} \!\cdot\! \frac{z_j^{(2)}}{p^{(d_2-1) / 2}}, \,
\textrm{\upshape and}~ \,
p\frac{p^{(d_1-1) / 2}}{z_i^{(1)}} \!\cdot\! \frac{p^{(d_2-1) / 2}}{z_j^{(2)}} \, .$$
It~shows that the product
$\smash{\prod_{i,j} (-1 - z_i^{(1)} \!z_j^{(2)} / p^{\frac{d_1 + d_2}2})}$
is
$p^{N_1N_2/4}$
times a~square. Indeed,~the sums occurring on the right hand side of~(\ref{quadm}) form a
$\Gal(\overline\bbQ/\bbQ)$-invariant
set of
$N_1N_2/4$~elements.
}
\eop
\end{lemma}

\begin{proposition}
\label{odd}
Let\/~$X$
be a smooth proper variety of odd
dimension\/~$d$
over a finite
field\/~$\bbF_{\!p^k}$
for\/~$p$
a prime and\/
$k$~odd.
Suppose~that\/
$p > 1 + \dim H^d_\et(X_{\overline\bbF_{\!p}}\!, \bbQ_l)$,
$$\smash{\dim H^d_\et(X_{\overline\bbF_{\!p}}\!, \bbQ_l) \equiv 2 \pmod 4 \, ,}$$
and that all eigenvalues
of\/~$\Frob$
on the cohomology are of\/
\mbox{$p$-adic}
valuation\/~$dk/2$.
I.e.,~that the Newton polygon has constant
slope\/~$d/2$.\smallskip

\noindent
Then\/~$\pm\sqrt{-p^{dk}}$
are among the~eigenvalues.\medskip

\noindent
{\bf Proof.}
{\em
Let~$C$
be the base extension of a supersingular elliptic curve defined
over~$\bbF_{\!p}$.
Such~do exist by the work of M.\,Eichler~\cite{Ei}, see also~\cite[Proposition~2.4, together with (1.10) and~(1.11)]{Gr}. The~eigenvalues
of~$\Frob$
on
$\smash{H^1_\et(C_{\overline\bbF_{\!p}}\!,\bbQ_l)}$
are~$\smash{\pm\sqrt{-p^k}}$.
Indeed,~the assumptions imply
$p \geq 5$,
and hence
$\smash{p > 2\sqrt{\mathstrut p}}$.
The~claim follows from Hasse's~bound.

For~$\Phi$
the characteristic polynomial
of~$\Frob$
on
$$V := [H^d_\et(X_{\overline\bbF_{\!p}}\!,\bbQ_l) \otimes H^1_\et(C_{\overline\bbF_{\!p}}\!,\bbQ_l)](\textstyle\frac{d+1}2) \, ,$$
Lemma~\ref{comb} guarantees that
$p \Phi(-1)$
is a perfect~square. But~all eigenvalues
of~$\Frob$
on~$V$
are
\mbox{$p$-adic}
units. As~they are
\mbox{$l$-adic}
units for every prime
$l \neq p$,
too, they must be roots of unity~\cite[Sec.\,18, Lemma~2]{Ca2}.

To~prove the assertion, we need to show
$\Phi(-1) = 0$.
Assuming~the contrary, we see from Lemma~\ref{roots} that, for some
$e \geq 1$,
all primitive
\mbox{$2p^e$-th}
roots of unity must be eigenvalues
of~$\Frob$
on~$V$.
As~with~$z$,
$(-z)$
is an eigenvalue, too, this enforces
$\dim V \geq 2(p-1)$,
that
is~$\dim H^d_\et(X_{\overline\bbF_{\!p}}\!,\bbQ_l) \geq p-1$,
a~contradiction.
}
\eop
\end{proposition}

\begin{remark}
In~principle, the idea behind this proof is to apply Theorem~\ref{main2} to
$X \times C$.
This~is, however, not sufficient as there may be
eigenvalues~$(-1)$
on the products
$[H^{i_1}(X,\bbQ_l) \otimes H^{i_2}(C,\bbQ_l)](\frac{d+1}2)$
for
$i_1 + i_2 = d+1$,
$i_2 \neq 1$.
\end{remark}

\begin{example}
\label{Stephan}
Proposition~\ref{odd} may fail in small characteristic, as is seen from the elliptic curve
$C$
over~$\bbF_{\!3}$,
given by
$y^2 = x^3 - x - 1$.
Then~$\#C(\bbF_{\!3}) = 1$,
which shows that
$C$
is supersingular. Moreover,~the characteristic polynomial of
$\Frob$
is
$T^2 - 3T + 3$.
The~eigenvalues are the two primitive twelfth roots of unity with positive real part, multiplied
by~$\sqrt{3}$.
\end{example}

\begin{remark}
\label{Stephan2}
When trying to carry over the argument from the proof, it turns out that,~on
$\smash{V := [H^1_\et(C_{\overline\bbF_{\!3}}\!,\bbQ_l) \times H^1_\et(C_{\overline\bbF_{\!3}}\!,\bbQ_l)](1)}$,
there are the primitive sixth roots of unity occurring as eigenvalues, together with
$1$,
which appears~twice. Therefore,
$\Phi(T) = (T^2 - T + 1)(T-1)^2$
and~$\Phi(-1) = 12$.

Alternatively, one might combine
$C$
with
$C'\colon y^2 = x^3 - x$,
which is supersingular having four~points. On the product
$\smash{[H^1_\et(C_{\overline\bbF_{\!3}}\!,\bbQ_l) \otimes H^1_\et(C'_{\overline\bbF_{\!3}}\!,\bbQ_l)](1)}$,
one finds
$\Phi(T) = (T^2 - T + 1)(T^2 + T + 1)$
and~$\Phi(-1) = 3$.
\end{remark}

\begin{lemma}
\label{roots}
Let\/~$\Phi \in \bbQ[T]$
be a monic polynomial such that all its roots are roots of~unity. Assume~that\/
$|\Phi(-1)| \neq 0$
is a multiple of a prime\
number\/~$p$.\smallskip

\noindent
Then\/~$\Phi$
is divisible by the cyclotomic polynomial\/
$\phi_{2p^e}$,
for
some\/~$e \geq 1$.
In~particular, if\/
$p=2$
then\/
$\Phi$
is divisible by\/
$\phi_{2^e}$,
for some\/
$e \geq 2$.\medskip

\noindent
{\bf Proof.}
{\em
$\Phi$~is
a product of cyclotomic
polynomials~$\phi_n$.
For~these, it is well known~\cite[Section~3]{Mo} that
$\phi_1(-1) = -2$,
$\phi_{2}(-1) = 0$,
$\phi_{2p^e}(-1) = p$
for
$p$
any prime number and
$e \geq 1$,
and
$\phi_n(-1) = 1$
in all other~cases. The~assertion follows directly from~this.
}
\eop
\end{lemma}

\subsection{The even-dimensional supersingular~case}

\begin{proposition}
\label{m1}
Let\/~$X$
be a smooth proper variety of even
dimension\/~$d$
over a finite
field\/~$\bbF_{\!p^k}$
for a
prime\/~$p \neq 2$
and\/
$k$~odd.
Suppose~that\/
$a(X) \equiv 1 \pmod 2$
and that all eigenvalues
of\/~$\Frob$
on\/~$H^d_\et(X_{\overline\bbF_{\!p}}\!,\bbQ_l(d/2))$
are\/
\mbox{$p$-adic}~units.
I.e.,~that the Newton polygon has constant
slope\/~$d/2$.

\begin{abc}
\item
Then\/~$(-1)$
is an eigenvalue or, for some\/
$e \geq 1$,
the primitive\/
$2p^e$-th
roots of unity are~eigenvalues.
\item
If\/~$p > \dim H^d_\et(X_{\overline\bbF_{\!p}}\!,\bbQ_l(d/2)) + 1$
or\/
$p > \dim H^d_\et(X_{\overline\bbF_{\!p}}\!,\bbQ_l(d/2))$
and\/
$X$
is projective then\/
$(-1)$
is an~eigenvalue.
\end{abc}\medskip

\noindent
{\bf Proof.}
{\em
The~eigenvalues are
\mbox{$l$-adic}
units, too, for every prime number
$l \neq p$,
hence they are roots of~unity. Theorem~\ref{main2} ensures that
$p \Phi(-1)$
is a perfect~square. Applying Lemma~\ref{roots} again, we immediately obtain assertion~a).

For~b), assume the~contrary. Then,~as eigenvalues, we have the
$(p-1)p^{e-1} \geq p-1$
primitive
$2p^e$-th
roots of unity and, in the projective case, the
number~$1$.
Thus,~altogether, there are at
least
$p-1$,
respectively
$p$,
of~them.
}
\eop
\end{proposition}

\begin{corollary}
\label{m1k3}
Let\/~$X$
be a supersingular\/
$K3$~surface
over a finite
field\/~$\bbF_{\!p^k}$
for\/
$p > 19$
a prime and\/
$k$~odd.
Then\/~$(-1)$
is an eigenvalue
of\/~$\Frob$
on\/
$H^2_\et(X_{\overline\bbF_{\!p}}\!,\bbQ_l(1))$.\medskip

\noindent
{\bf Proof.}
{\em
For~$K3$~surfaces,
the Hodge spectral sequence degenerates
at~$E_1$~\cite[Proposition~1.1.a)]{De3}
and, hence, the conjugate spectral sequence degenerates
at~$E_2$~\cite[Lemma~8.27.2]{BO}.
Moreover,~all
$H^i(X/W)$
are torsion-free (\cite[II.7.2]{Il} or~\cite[Proposition~1.1.c)]{De3}).
Consequently,~we
have~$a(X) = \dim H^2(X,\calO_X) = 1$.
Cf.~Remark~\ref{techn}.b). The~claim now follows from Proposition~\ref{m1}.b).
}
\eop
\end{corollary}

\begin{rems}
\label{Artin}
\begin{iii}
\item
Corollary~\ref{m1k3} refines the observation of M.\,Artin \mbox{\cite[6.8]{Ar}} that the field of definition of the
\mbox{rank-$22$}
Picard group always
contains~$\bbF_{\!p^2}$.
\item
More generally, let
$X$
be any
$K3$~surface
over a finite field
$\bbF_{\!p^k}$
for\/~$p \neq 2$
a prime. Theorem~\ref{main2} then asserts that,
for~$z_1, \ldots, z_{22}$
the eigenvalues
of~$\Frob$
on
$\smash{H^2_\et(X_{\overline\bbF_{\!p}}\!, \bbQ_l(1))}$,
the expression
$p^k \Phi(-1) = p^k (1+z_1)\ldots(1+z_{22})$
is always a square
in~$\bbQ$.
\item
Corollary~\ref{m1k3} is clearly false in small characteristic, as may be seen from the example~below.
\end{iii}
\end{rems}

\begin{exs}
\begin{iii}
\item
For~$C$
the supersingular elliptic curve
over~$\bbF_{\!3}$
from Example~\ref{Stephan}, put
$X := \Kum(C \times C)$.
Then
$$\Phi(T) = (T^2-T+1)(T^2+T+1)^5(T-1)^{10}.$$
In~particular, we have
$\Phi(-1) = 3 \!\cdot\! 2^{10}$,
in agreement with Remark~\ref{Artin}.ii).

Indeed,~we saw that the characteristic polynomial of
$\Frob$
on
$\smash{H^2_\et((C \times C)_{\overline\bbF_{\!3}}\!, \bbQ_l(1))}$
is~$(T^2-T+1)(T-1)^4$.
Furthermore,~the 16 two-torsion points form five orbits of size three together with the~origin.
\item
Let~$X := \Kum(C \times C')$
be the Kummer surface associated to the product of the two supersingular elliptic curves
over~$\bbF_{\!3}$,
considered in Remark~\ref{Stephan2}. Then,~once again,
$\Phi(T) = (T^2-T+1)(T^2+T+1)^5(T-1)^{10}$.

In~fact, we saw that the characteristic polynomial of
$\Frob$
on
$\smash{H^2_\et((C \times C')_{\overline\bbF_{\!3}}\!, \bbQ_l(1))}$
is
$(T^2-T+1)(T^2+T+1)(T-1)^2$.
In~addition, four two-torsion points are defined over the base field, while the others form four orbits of size~three.
\end{iii}
\end{exs}

\subsection{Surfaces. The Artin-Tate formula}

For~surfaces, the assertion of Theorem~\ref{main} is implied by the Tate~conjecture. More~precisely,

\begin{proposition}
\label{ArT}
Let\/~$X$
be a smooth projective surface over a finite
field\/~$\bbF_{\!q}$
of
characteristic\/~$p$
and
let\/~$\smash{\Phi = \Phi_1^{(2)} \in \bbQ[T]}$
be the characteristic polynomial
of\/~$\Frob$
on\/
$\smash{H^2_\et(X_{\overline\bbF_{\!q}}\!, \bbQ_l(1))}$.
Put\/~$N := \deg \Phi$~and
$$\alpha(X) := \dim H^2(X,\calO_X) - \dim H^1(X,\calO_X) + \textstyle\frac12 \dim H^1_\et(X_{\overline\bbF_{\!q}}\!, \bbQ_l) \, .$$
Suppose~that the Tate~conjecture is true
for\/~$X$.\smallskip

\noindent
Then\/
$(-2)^N q^{\alpha(X)} \Phi(-1)$
is a square
in\/~$\bbQ$.\medskip

\noindent
{\bf Proof.}
{\em
Denote~the zeroes
of~$\Phi$,
i.e.~the eigenvalues of
$\Frob$,
by
$z_1,\ldots,z_\varrho = 1$,
$z_{\varrho+1},\ldots,z_N \neq 1$.
If~$\Phi(-1) = 0$
then the assertion is true,~trivially. Thus,~let us suppose the contrary from now~on. Then~the zeroes
$z_i \neq 1$
come in pairs of complex conjugate numbers. In~particular,
$N-\varrho$
is~even.

Furthermore,~$\Frob$
and~$\Frob^2$
have the
eigenvalue~$1$
with the same~multiplicity. Hence, the Tate conjecture predicts the rank
of~$\smash{\Pic(X_{\bbF_{\!q^2}})}$
not to be higher than that
of~$\Pic(X)$.
This~shows that
$\smash{X_{\bbF_{\!q^2}}}$,
too, fulfills the Tate~conjecture.

We~are therefore in a situation where the Artin-Tate formula~\cite[Theorem~6.1]{Mi} computes the discriminants of the Picard lattices
$\Pic(X)$
and
$\smash{\Pic(X_{\bbF_{\!q^2}})}$,
at least up to square~factors. The~results are
$$(-1)^{\varrho-1} q^{\alpha(X)} \!\!\!\prod_{i=\varrho+1}^N\!\! (1-z_i)
\quad {\rm and} \quad
(-1)^{\varrho-1} q^{2\alpha(X)} \!\!\!\prod_{i=\varrho+1}^N\!\! (1-z_i^2) \, .$$
Moreover,~equality of the ranks implies that
$\smash{\disc \Pic(X) / \disc \Pic(X_{\bbF_{\!q^2}})}$
is a necessarily perfect~square. This~is a standard observation from the theory of~lattices. We~conclude that
$q^{\alpha(X)} \prod_{i=\varrho+1}^N (1+z_i)$
is a square
in~$\bbQ$.

On~the other hand,
$\smash{(-2)^N q^{\alpha(X)} \Phi(-1) = 2^{N+\varrho} q^{\alpha(X)} \prod_{i=\varrho+1}^N (1+z_i)}$
such that the assertion follows from the fact that
$\varrho \equiv N \pmod 2$.
}
\eop
\end{proposition}

\begin{rems}
\begin{iii}
\item
The~Artin-Tate formula appears to us as a very natural consequence of the Tate conjecture and the cohomological~machinery. Thus,~we find it very astonishing that it has the potential to produce incompatible results for a variety and its base~extension.

Of~course, this does not happen for polynomials that really occur as the characteristic polynomial of the Frobenius on a certain~variety. But~it occurs for polynomials that otherwise look~plausible. This~observation was actually the starting point of our~investigations.
\item
One~might want to compare the Picard lattice
of~$X$
with that of
$X_{q^n}$
for~$n > 2$.
But~this leads to nothing~new~\cite[Corollary~18.i)]{EJ1}.
\item
Suppose~$\smash{\dim H^1(X,\calO_X) = \frac12 \dim H^1_\et(X_{\overline\bbF_{\!q}}\!, \bbQ_l)}$
and that
$X$
fulfills the assumptions of Remark~\ref{techn}.b).
Then~$\smash{\alpha(X) = a(X) = \dim H^2(X,\calO_X)}$.
We~do not know how closely Milne's invariant
$\alpha(X)$
and our
invariant~$a(X)$
are related for ``pathological''~surfaces.
\end{iii}
\end{rems}

\begin{example}
\label{k3}
Let~$X$
be the double cover
of~$\bP^2_{\!\bbF_{\!7}}$,
given~by
$$w^2 = 6x^6 + 6x^5y + 2x^5z + 6x^4y^2 + 5x^4z^2 + 5x^3y^3 + 
x^2y^4 + 6xy^5 + 5xz^5 + 3y^6 + 5z^6 \, .$$
This~is a
$K3$~surface
of degree~two.

The~numbers of points
on~$X$
over the finite fields
$\bbF_{\!7}, \ldots, \bbF_{\!7^{10}}$
are
$60$,
$2\,488$,
$118\,587$,
$5\,765\,828$,
$282\,498\,600$,
$13\,841\,656\,159$,
$678\,225\,676\,496$, 
$33\,232\,936\,342\,644$,
$1\,628\,413\,665\,268\,026$,
and
$79\,792\,266\,679\,604\,918$.

For~the characteristic polynomial
of~$\Frob$
on~$\smash{H^2_\et(X_{\overline\bbF_{\!7}}\!, \bbQ_l(1))}$,
this information leaves us with two candidates, one for each sign in the functional~equation,
\begin{align*}
\Phi_i(t) &=
\frac17
\big(
7 t^{22} - 10 t^{21} + t^{20} - t^{19} + 6 t^{18} - 3 t^{17} - 2 t^{16} + 4 t^{14} - t^{13} - t^{12}\\[-2mm]
& \hspace{1.8cm} ~{} + (-1)^i
(- t^{10} - t^9 + 4 t^8 - 2 t^6 - 3 t^5 + 6 t^4 - t^3 + t^2 - 10 t + 7)
\big)
\end{align*}
\noindent
for~$i = 0, 1$.
All roots are of absolute
value~$1$.

However,~$\Phi_0(-1) = 60/7$
and~$(-2)^N 7^{a(X)} \Phi_0(-1) = 2^{24} \!\cdot\! 3 \!\cdot\! 5$
is a non-square, which contradicts Theorem~\ref{main}.
Therefore,~$\Phi_1$
is the characteristic polynomial
of
$\Frob$
on
$\smash{H^2_\et(X_{\overline\bbF_{\!7}}\!, \bbQ_l(1))}$.
The~minus sign holds in the functional~equation.
\end{example}

\begin{remark}
Alternatively, we may argue as follows.
Assume~$\Phi_0$
is the characteristic polynomial. Then
$\smash{\rk \Pic(X_{\bbF_{\!49}}) = \rk \Pic(X) = 2}$.
Indeed,~the Tate conjecture is proven for
$K3$~surfaces
in characteristic
$\geq\!3$~\cite[Corollary~2]{Char}, \cite[Theorem~1]{Pe}, cf.~\cite{LMS}. The~Artin-Tate formula states that
$\smash{\disc \Pic(X_{\bbF_{\!49}}) \in (-465) (\bbQ^*)^2}$
and
$\disc \Pic(X) \in (-31) (\bbQ^*)^2$.
As
$\smash{\frac{-465}{-31} = 15}$
is a non-square, this is~contradictory.
\end{remark}

\subsection{Cubic fourfolds}

\begin{example}
\label{cub4f}
Let\/~$X$
be the subvariety
of\/~$\bP^5_{\!\bbF_{\!2}}$,
given~by
\begin{align*}
 & x_0^3 + x_0^2x_1 + x_0^2x_4 + x_0^2x_5 + x_0x_1x_2 + x_0x_1x_3 + x_0x_1x_4 + x_0x_2x_3 + x_0x_2x_4 \\
 & {} \hspace{-0.1cm} + x_0x_3x_4 + x_0x_4^2 + x_0x_4x_5 + x_0x_5^2 + x_1^3 + x_1x_2^2 + x_1x_2x_4 + x_1x_2x_5 + x_1x_3^2 \\
 & {} \hspace{0.7cm} + x_1x_3x_5 + x_1x_4^2+ x_1x_4x_5 + x_2^3 + x_2^2x_5 + x_2x_3^2 + x_3^2x_4 + x_3^2x_5 + x_3x_4^2 \\
 & {} \hspace{6.8cm} + x_3x_5^2 + x_4^3 + x_4^2x_5 + x_4x_5^2 + x_5^3 = 0 \, .
\end{align*}
This~is a smooth cubic~fourfold. We~have
$\dim H^4(X,\calO_X) = 0$,
$\dim H^3(X,\Omega_X^1) = 1$,
and
$\dim H^2(X,\Omega_X^2) = 21$,
According~to Remark~\ref{techn}.b), this shows
$N = 23$
and~$a(X) = 1$.

The~numbers of points
on~$X$
over the finite fields
$\bbF_{\!2}, \ldots, \bbF_{\!2^{11}}$
are
$33$,
$361$,
$4\,545$,
$69\,665$,
$1\,084\,673$,
$17\,044\,609$,
$270\,543\,873$,
$4\,311\,990\,785$,
$68\,853\,026\,817$,
$1\,100\,586\,076\,161$,
and
$17\,600\,769\,409\,025$.
The~characteristic polynomial
of~$\Frob$
on
$\smash{H^4_\et(X_{\overline\bbF_{\!2}}\!, \bbQ_l(2))}$
is
\begin{align*}
\Phi(t) &= \frac12 (t-1) (2t^{22} - t^{21} - t^{20} + 2t^{19} - 2t^{17} + t^{16} + t^{15} - 2t^{14} + t^{13} + t^{12} \\[-2mm]
& \hspace{2.8cm} ~{} - t^{11} + t^{10} + t^9 - 2t^8 + t^7 + t^6 - 2t^5 + 2t^3 - t^2 - t + 2) \, . 
\end{align*}
It~turns out
that~$\Phi(-1) = -1$,
in agreement with Theorem~\ref{main}. Observe~that, in this example, the assertion of Theorem~\ref{main2} is true, albeit the characteristic of the base field
is~$2$.
\end{example}

\begin{remark}
The~degree
$22$~factor
of~$\Phi$
is irreducible
over~$\bbQ$.
In~particular,
$X$
is certainly not special in the sense of B.~Hassett~\cite{Ha}.
\end{remark}

\parindent0pt

\begin{thebibliography}{99}

\bibitem{Ar}
\textit{M.~Artin},
Supersingular
$K3$
surfaces, Ann.\ Sci.\ \'Ecole Norm.\ Sup.\ \textbf{7} (1974), 543--567.

\bibitem{SGA4}
\textit{M.~Artin, A.~Grothendieck,} et \textit{J.-L.~Verdier} (avec la collaboration de \textit{P.~Deligne} et \textit{B.~Saint-Donat}),
Th\'eorie des topos et cohomologie \'etale des sch\'emas, S\'eminaire de G\'eom\'etrie Alg\'ebrique du Bois Marie 1963--1964 (SGA\,4), Lecture Notes in Math.~269, 270, 305, Springer, Berlin, Heidelberg, New York~1972--1973.

\bibitem{Bea}
\textit{A.~Beauville},
Complex algebraic surfaces, LMS Lecture Note Series~68, Cambridge University Press, Cambridge~1983.

\bibitem{Be}
\textit{P.~Berthelot},
Cohomologie cristalline des sch\'emas de caract\'eristique
$p > 0$,
Lecture Notes in Math.~407, Springer, Berlin, New York~1974.

\bibitem{BO}
\textit{P.~Berthelot} and \textit{A.~Ogus},
Notes on crystalline cohomology, Princeton University Press, Princeton~1978.

\bibitem{Ca1}
\textit{J.\,W.\,S.~Cassels},
Arithmetic on curves of genus~1. VIII. On the conjectures of Birch and Swinnerton-Dyer, J.~f\"ur die Reine und Angew.\ Math.\ \textbf{217} (1965), 180--189.

\bibitem{Ca2}
\textit{J.\,W.\,S.~Cassels},
Global fields, in: Algebraic number theory, Edited by J.~W.~S.~Cassels and A.~Fr\"ohlich, Academic Press and Thompson Book Co., London and Washington 1967, 42--84.

\bibitem{Cham}
\textit{A.~Chambert-Loir},
Cohomologie cristalline: un survol, Exposition.\ Math.\ \textbf{16} (1998), 333--382.

\bibitem{Char}
\textit{F.~Charles}, The Tate conjecture for
$K3$~surfaces
over finite fields, {\tt arXiv:1206.4002}.

\bibitem{SGA41/2}
\textit{P.~Deligne} (avec la collaboration de \textit{J.~F.~Boutot}, \textit{A.~Grothendieck, L.~Illusie}, et \textit{J.-L.~Verdier}),
Cohomologie \smash{\'Etale,} S\'eminaire de G\'eom\'etrie Alg\'ebrique du Bois Marie (SGA\,4$\frac{1}{2}$), Lecture Notes in Math.~569, Springer, Berlin, Heidelberg, New York~1977.

\bibitem{De1}
\textit{P.~Deligne},
La conjecture de Weil~I, Publ.\ Math.\ IHES \textbf{43} (1974), 273--307.

\bibitem{De2}
\textit{P.~Deligne},
La conjecture de Weil~II, Publ.\ Math.\ IHES \textbf{52} (1980), 137--252.

\bibitem{De3}
\textit{P.~Deligne},
Rel\`evement des surfaces
$K3$
en caract\'eristique nulle, Prepared for publication by Luc Illusie, in: Algebraic surfaces (Orsay 1976--78), Lecture Notes in Math.~868, Springer, Berlin-New York~1981, 58--79.

\bibitem{DI}
\textit{P.~Deligne} and \textit{L.~Illusie},
Rel\`evements modulo
$p^2$
et d\'ecomposition du complexe de de Rham, Invent.\ Math.\ \textbf{89} (1987), 247--270.

\bibitem{SGA7}
\textit{P.~Deligne} and \textit{N.~Katz},
Groupes de Monodromie en G\'eom\'etrie Alg\'ebrique, S\'e\-mi\-naire de G\'eom\'etrie Alg\'ebrique du Bois Marie 1967--1969 (SGA\,7), Lecture Notes in Math.~288, 340, Springer, Berlin, Heidelberg, New York~1973.

\bibitem{Di}
\textit{J.~Dieudonn\'e},
\'El\'ements d'analyse, Tome~II, 3${}^{\rm e}$ \'edition, Gauthier-Villars, Paris~1983.

\bibitem{Ei}
\textit{M.~Eichler},
Zur Zahlentheorie der Quaternionen-Algebren, J.~f\"ur die Reine und Angew.\ Math.\ \textbf{195} (1955), 127--151.

\bibitem{EJ1}
\textit{A.-S.~Elsenhans} and \textit{J.~Jahnel},
On Weil polynomials of
$K3$~surfaces,
in: Algorithmic Number Theory (ANTS~9), Lecture Notes in Computer Science~6197, Springer, Berlin~2010, 126--141.

\bibitem{EJ2}
\textit{A.-S.~Elsenhans} and \textit{J.~Jahnel},
The Picard group of a
$K3$~surface
and its reduction
modulo~$p$,
Algebra \& Number Theory \textbf{5} (2011), 1027--1040.

\bibitem{EJ3}
\textit{A.-S.~Elsenhans} and \textit{J.~Jahnel},
On the computation of the Picard group for
$K3$~surfaces,
Mathematical Proceedings of the Cambridge Philosophical Society \textbf{151} (2011), 263--270.

\bibitem{Gr}
\textit{B.\,H.~Gross},
Heights and the special values of
\mbox{$L$-series},
in: Number theory (Montreal 1985), CMS Conf.\ Proc.\,7, AMS, Providence~1987, 115--187.

\bibitem{Ha}
\textit{B.~Hassett},
Special cubic fourfolds, Compositio Math.\ \textbf{120} (2000), 1--23.

\bibitem{Ho}
\textit{T.~Honda},
Isogeny classes of abelian varieties over finite fields, J.\ Math.\ Soc.\ Japan \textbf{20} (1968), 83--95.

\bibitem{Il}
\textit{L.~Illusie},
Complexe de de Rham-Witt et cohomologie cristalline, Ann.\ Sci.\ \'Ecole Norm.\ Sup.\ \textbf{12} (1979), 501--661.

\bibitem{IR}
\textit{L.~Illusie} and \textit{M.~Raynaud},
Les suites spectrales associ\'ees au complexe de de Rham-Witt, Publ.\ Math.\ IHES \textbf{57} (1983), 73--212.

\bibitem{KM}
\textit{N.\,M.~Katz} and \textit{W.~Messing},
Some consequences of the Riemann hypothesis for varieties over finite fields, Invent.\ Math.\ \textbf{23} (1974), 73--77.

\bibitem{LMS}
\textit{M.~Lieblich, D.~Maulik}, and \textit{A.~Snowden}, Finiteness of
$K3$~surfaces
and the Tate conjecture, {\tt arXiv:\discretionary{}{}{}1107.1221}.

\bibitem{vL}
\textit{R.~van Luijk},
$K3$~surfaces
with Picard number one and infinitely many rational points, Algebra \& Number Theory, \textbf{1} (2007), 1--15.

\bibitem{Man}
\textit{Yu.\,I.~Manin},
Theory of commutative formal groups over fields of finite characteristic (Russian), Uspehi Mat.\ Nauk \textbf{18} (1963), 3--90.

\bibitem{Maz1}
\textit{B.~Mazur},
Frobenius and the Hodge filtration, Bull.\ Amer.\ Math.\ Soc.\ \textbf{78} (1972), 653--667.

\bibitem{Maz2}
\textit{B.~Mazur},
Frobenius and the Hodge filtration (estimates), Ann.\ of Math.\ \textbf{98} (1973), 58--95.

\bibitem{Mi}
\textit{J.\,S.~Milne},
On a conjecture of Artin and Tate, Ann.\ of Math.\ \textbf{102} (1975), 517--533.

\bibitem{MS}
\textit{J.\,W.~Milnor} and \textit{J.\,D.~Stasheff},
Characteristic classes, Annals of Mathematics Studies~76, Princeton University Press, Princeton~1974.

\bibitem{Mo}
\textit{K.~Motose},
On values of cyclotomic polynomials~VIII, Bull.\ Fac.\ Sci.\ Technol.\ Hirosaki Univ.\ \textbf{9} (2006), 15--27.

\bibitem{Pe}
\textit{K.\,P.~Pera},
The Tate conjecture for
$K3$~surfaces
in odd characteristic, {\tt arXiv:1301.6326}.

\bibitem{Se}
\textit{J.-P.~Serre},
Sur la topologie des vari\'et\'es alg\'ebriques en
caract\'eristique~$p$,
in: Symposium internacional de topolog\'\i a algebraica, Universidad Nacional Aut\'onoma de M\'exico and UNESCO, 24--53.

\bibitem{Si}
\textit{J.\,H.~Silverman},
The arithmetic of elliptic curves, Graduate Texts in Mathematics~106, Springer, New York~1986.

\bibitem{Sp}
\textit{E.\,H.~Spanier},
Algebraic topology, McGraw-Hill, New York, To\-ron\-to, Lon\-don~1966.

\bibitem{St}
\textit{U.~Stuhler},
Konjugationsklassen unipotenter Elemente in einfachen algebraischen Grup\-pen vom Typ
$B_n$,
$C_n$,
$D_n$
und~$G_2$,
Dissertation, G\"ot\-tin\-gen~1970.

\bibitem{Su}
\textit{J.~Suh},
Symmetry and parity in Frobenius action on cohomology, Compositio Math.\ \textbf{148} (2012), 295--303.

\bibitem{Ur}
\textit{T.~Urabe},
The bilinear form of the Brauer group of a surface, Invent.\ Math.\ \textbf{125} (1996), 557--585.

\bibitem{Wa}
\textit{G.\,E.~Wall},
On the conjugacy classes in the unitary, symplectic and orthogonal groups, J.\ Austral.\ Math.\ Soc.\ \textbf{3} (1963), 1--62.

\bibitem{Zar}
\textit{Yu.\,G.~Zarhin},
The Brauer group of an abelian variety over a finite field (Russian), Izv.\ Akad.\ Nauk SSSR Ser.\ Mat.\ \textbf{46} (1982), 211--243.

\bibitem{Zag}
\textit{D.\,B.~Zagier},
Zetafunktionen und quadratische K\"orper, Springer, Berlin, New York~1981.

\end{thebibliography}
\end{document}